\documentclass{article}

\title{A simple proof of a duality theorem with applications in viscoelasticity}
\date{}
\author{Andrzej Hanyga \\
ul. Bitwy Warszawskiej 1920 r 14/52\\
02-366 Warszawa, Poland}

\usepackage{hanmin,amsmath,amssymb}
\newtheorem{theorem}{Theorem}
\newcommand{\D}{\mathrm{D}}

\begin{document}

\maketitle

\section*{Notations}
$[0,\infty[ := \{ x \in \mathbf{R} \mid 0 \leq x < \infty \}$\\
$]0,\infty[ := \{ x \in \mathbf{R} \mid 0 < x < \infty \}$

\begin{abstract}
A new concise proof is given of a duality theorem connecting completely monotone relaxation functions 
with Bernstein class creep functions. The proof makes use of the theory of complete Bernstein functions and 
Stieltjes functions and is based on a relation between these two function classes. 
\end{abstract}

\textbf{Keywords:}
{viscoelasticity,completely monotone, Bernstein, complete Bernstein, Stieltjes}

\section{Introduction.}

In a few papers Hanyga and Seredy\'{n}ska \cite{SerHanJMP,HanWM} demonstrated the utility of 
the concept of complete Bernstein functions (CBFs) in the theory of wave 
propagation and attenuation in viscoelastic media with completely monotone relaxation moduli. It is 
therefore reasonable to apply complete Bernstein functions 
in other contexts of the theory of linear viscoelasticity as well. 

In this note I shall 
demonstrate the application of the CBFs and their inseparable fellow travellers the Stieltjes functions 
in the proof of a well known relation between locally integrable completely monotone (LICM) 
relaxation moduli $R$ and and Bernstein class creep functions $C$. I recall the viscoelastic
duality relation \cite{Molinari73,HanISIMM,MainardiBook}
\begin{equation} \label{duality}
\tilde{R}(p) \, \tilde{C}(p) = p^{-2}
\end{equation}
where 
\begin{equation}\label{Laplace}
\tilde{f}(p) := \int_0^\infty \e^{-p t} \, f(t) \, \dd t
\end{equation}
denotes the Laplace transform
for any locally integrable function $f$ such that the transform exists.

It was previously shown by other methods \cite{Molinari73,HanISIMM} that, apart from some singular terms, 
if $R$ is LICM  then $C$ is a Bernstein function and conversely.
The new proof is more elegant and concise than the 
previously given proofs \cite{Molinari73,HanISIMM} of the same relation. 

\section{Preliminaries.}

It is convenient for our considerations to consider functions $f$ in 
$\mathcal{L}^1_{\mathrm{loc}}([0,\infty[)$ as convolution operators $f\ast$ mapping 
$g \in \mathcal{L}^1_{\mathrm{loc}}([0,\infty[)$ to $f\ast g \in \mathcal{L}^1_{\mathrm{loc}}[0,\infty[)$.

The convolution of two locally integrable functions $f$ and $g$ on $[0,\infty[$
is defined by the formula
\begin{equation}
(f\ast g)(t) := \int_0^t f(s) \, g(t-s) \, \dd s
\end{equation}
If the Laplace transforms $\tilde{f}(p)$ and $\tilde{g}(p)$ exist for some $p \geq 0$, then the convolution 
$f\ast g$ also has the Laplace transform at $p$ and 
\begin{equation} \label{LaplId}
(f\ast g)\,\widetilde{}\,(p) = \tilde{f}(p)\, \tilde{g}(p) 
\end{equation}

We shall also need an identity operator on $\mathcal{L}^1_{\mathrm{loc}}([0,\infty[)$:
\begin{equation} \label{unity}
\mathrm{U}\, f = f
\end{equation}
For the sake of convenience we shall also write \eqref{unity} in the form
\begin{equation}
u\ast f = f
\end{equation}

Extending \eqref{LaplId} to \eqref{unity} we have
$\tilde{u}(p) \tilde{f}(p) = \tilde{f}(p)$, whence
\begin{equation}
\tilde{u}(p) = 1,\; p \geq 0
\end{equation}

Note that if $\sigma = R\ast \dot{\epsilon}$, where $\sigma$ represents the stress, 
$\dot{\epsilon}$ is the strain rate and the relaxation modulus $R = N \, u + f_0$, where 
$f_0 \in \mathcal{L}^1_{\mathrm{loc}}([0,\infty[)$, then $\sigma = N \, \dot{\epsilon} + 
f_0\ast \dot{\epsilon}$. In the absence of the second term the first term  represents Newtonian viscosity.
We shall however see that for the validity of the duality relation the appearance of a term i$ $b \, u is necessary.

The mathematical function classes required in the proof of the duality theorem are 
explained in the appendix along with their properties which will be needed in the
following.

\section{The proofs.}

\begin{theorem}
If $f_0$ is LICM and not identically zero and 
\begin{equation} \label{f0}
f(t) = f_0(t) + \beta \, u(t)
\end{equation}
where $\beta \geq 0$, then 
\begin{equation} \label{y2} 
1/[p \tilde{f}(p)] = p\, \tilde{h}(p)
\end{equation} 
where $h$ is a Bernstein function.

We also have $0 < f_0(0) \leq \infty$. If $\beta > 0$ or $f_0(t)$ is unbounded at 0 then
$h(0) = 0$, otherwise $h(0) = 1/f_0(0)$.

The limit $f_\infty : = \lim_{t\rightarrow \infty} f_0(t)$ exists and is non-negative. \\
If $f_\infty > 0$, then for $t \rightarrow \infty$ the function $h(t)$ tends to $1/f_\infty$, 
otherwise it diverges to infinity.
\end{theorem}

\noindent\textbf{Proof.} 

We shall use a few theorems on  of the CBFs and Stieltjes functions quoted in the appendix to construct the Bernstein 
function $h$.

Since $f_0$ is LICM, there is a Borel measure $\mu$ on $]0,\infty[$ satisfying the inequality
\begin{equation} \label{x2}
\int_{]0,\infty[} (1 + s)^{-1} \,\mu(\dd s) < \infty
\end{equation}
such that 
\begin{equation}
f_0(t) = \int_{]0,\infty[} \e^{-s t} \mu(\dd s)
\end{equation}
It follows that
\begin{equation}
\tilde{f}_0(p) = \int_{]0,\infty[} (s + p)^{-1} \, \mu(\dd s) 
\end{equation}
hence, by eq.~\eqref{app2}, $p\, \tilde{f}(p) = \beta\, p + p\, \tilde{f}_0(p)$ is a 
CBF and $1/[p \tilde{f}(p)]$  is a Stieltjes function.
Hence there are non-negative real numbers $a, b$ and a Borel measure $\nu$ on $]0,\infty[$ satisfying the ineqality
\begin{equation} \label{x0}
\int_{]0,\infty[} (1 + s)^{-1}\,\nu(\dd s) < \infty
\end{equation}
such that
\begin{equation}\label{x11}
1/[p \, \tilde{f}(p)] = a + \frac{b}{p} + \int_{]0,\infty[} \frac{\nu(\dd r)}{r + p} 
\end{equation}

The last term is the Laplace transform $\tilde{g}(p)$ of the LICM function
\begin{equation}
g(t) := \int_{]0,\infty[} \e^{-t r} \, \nu(\dd r)
\end{equation}

If
$$G(t) := \int_0^t g(s)\, \dd s + a$$
then
$$p \, \tilde{G}(p) = \tilde{g}(p) + a$$
and $G$ is a Bernstein function. 

Concerning the term $b/p$ appearing in \eqref{x11} we note that  
$b/p^2$ is the Laplace transform of the Bernstein function $b \, t, t \geq 0$. 
Hence eq.~\eqref{y2} holds with the Bernstein function 
\begin{equation} \label{ht}
h(t) := b\, t + G(t) 
\end{equation}

It remains to consider the limits of these functions.

The function $f_0(t) \geq 0$ is non-increasing and $f_0 \not\equiv 0$.  
$f_0(t)$ may be unbounded at 0, otherwise the limit $f_0(0)$ of $f_0(t)$ at 0 exists and $f_0(0) > 0$. On the other hand $h(0) \geq 0$
is always finite.
On account of equation~\eqref{y2} 
$$h(0)  = \lim_{p\rightarrow \infty} 1/[\beta \, p + p \, \widetilde{f_0}(p)] $$ 
If $\beta > 0$ or $f_0(t)$ is unbounded at 0, then $h(0) = 0$, otherwise $h(0) = 1/f_0(0)$.

The function $f_0$ satisfies the inequalities $0 \leq f_0(t) \leq f_0(1)$ for $t > 1$ and its non-increasing. 
Consequently the limit $f_\infty := \lim_{t\rightarrow \infty} f_0(t)$ exists and is non-negative. We now note that 
$$\lim_{p\rightarrow 0} [p \, \tilde{f}(p)] = \lim_{p\rightarrow 0} [p \, \widetilde{f_0}(p))] = f_\infty,$$
If $f_\infty > 0$ then 
$$\lim_{t\rightarrow \infty} h(t) = \lim_{p \rightarrow 0} [p \, \tilde{h}(p)] = 1/f_\infty$$
on account of \eqref{y2}. Otherwise $h(t)$ is unbounded at infinity and 
the limit $\lim_{t\rightarrow \infty} h(t)$ does not exist.
\\
$\mbox{ }$ \hfill $\Box$\\

\begin{theorem}
If $h$ is a BF not identically zero, then there is a LICM function $f_0$ and a real 
$b \geq 0$ such that 
\begin{equation} \label{y4}
1/[p\, \tilde{h}(p)] = p \, \tilde{f}(p)
\end{equation}
where
\begin{equation}
f(t) := b \, u(t) + f_0(t)
\end{equation}

The function $h(t)$ is either bounded and tends to a positive limit at infinity or it diverges to infinity.
In the first case we have the identity
$$\lim_{t\rightarrow \infty} f_0(t) = \lim_{t \rightarrow \infty} 1/h(t),$$
otherwise $f_0(t)$ tends to 0 at infinity.

If $h(0) > 0$, then $f_0(t)$ is bounded at 0,  $f_0(0) = 1/h(0)$ and $b = 0$,
otherwise either $b > 0$ or $f_0(t)$ is unbounded at 0. 

If $h(0) = 0$, then $b = 1/h^\prime(0)$ for $h^\prime(0) \geq 0$.

\end{theorem}

\noindent\textbf{Proof.}

Since $p \, \tilde{h}(p) = [h^\prime]\tilde(p) + h(0)$ and the derivative $h^\prime$ \
of $h$ is LICM, $p \, \tilde{h}(p)$ is a Stieltjes function and therefore 
$1/[p \, \tilde{h}(p)]$ is a CBF.

Eq.~\eqref{app2} implies that there are two reals $a, b \geq 0$ and a Borel measure $\nu$ satisfying 
\eqref{x0} such that
$$1/[p \, \tilde{h}(p)] = a + b\, p + p \int_{]0,\infty[} \frac{\nu(\dd r)}{r + p}$$
Let 
$$f_0(t) := a + \int_{]0,\infty[} \e^{-r t}\, \nu(\dd r)$$
and
$$f(t) := f_0(t) + b \, u(t)$$
$f_0$ is clearly LICM and eq.~\eqref{y4} is satisfied.

Furthermore, $\lim_{t\rightarrow \infty} h(t) > 0$ and 
\begin{equation}
\lim_{t \rightarrow \infty}\, 1/h(t) = \lim_{p \rightarrow 0} 1/[p \, \tilde{h}(p)]
= \lim_{p\rightarrow 0} [p \, \tilde{f}(p)] = \lim{p\rightarrow 0} \widetilde{f_0}(p) = \lim_{t\rightarrow \infty}\, f_0(t)
\end{equation}

At the other end we note that $h(0) \geq 0$ exists. If $b = 0$, and $f_0(t)$ is bounded at 0,  then
$$h(0) = \lim_{p\rightarrow\infty} [p \, \tilde{h}(p)] = \lim_{p\rightarrow\infty} 1/[b \, p + p\, \widetilde{f_0}(p)] 
= 1/f_0(0),$$ otherwise $h(0) = 0$.

Hence if $h(0) > 0$, then $b = 0$ and $f_0(t)$ tends to $1/h(0)$ for $t \rightarrow 0$,
while if $h(0) = 0$ then either $b > 0$ or $f_0(t)$ diverges to infinity at 0.

If $h(0) = 0$, then $\lim_{p\rightarrow \infty} [p^2 \, \tilde{h}(p)] = \lim_{p\rightarrow \infty} [p \, \widetilde{h^\prime}(p)]
 = \lim_{t\rightarrow 0} h^\prime(t)$. The last limit exists because $h^\prime$ is LICM, but it may be infinite. 
On the other hand \eqref{y2} implies that $\lim_{p\rightarrow \infty} [p^2 \, \tilde{h}(p)] =
\lim_{p\rightarrow \infty} 1/[b + \widetilde{f_0}(p)]$. We now note that
\begin{equation} \label{y1} 
\lim_{p\rightarrow \infty} \widetilde{f_0}(p) = \lim_{p\rightarrow \infty} \{p \, [1/p \widetilde{f_0}(p)]\} =
\lim_{t\rightarrow 0} \int_0^t f_0(t) \, \dd t = 0
\end{equation}

The last equation in \eqref{y1} follows from the fact that $f_0$ is integrable over $[0,1]$.
Indeed,  for $t \leq 1$  
$$\int_0^t f_0(t) \, \dd t = \int_0^1 f_0(s) \, \theta(t - s) \, \dd s$$
where $\theta$ denotes the unit step function, $f_0(s)  \, \theta(t - s) \leq f_0(s)$ and 
$f_0(s)  \, \theta(t - s) \rightarrow 0$ for $0 \leq s \leq 1$ as $t \rightarrow 0$,
hence the last equation in \eqref{y1} follows from the Lebesgue dominated convergence theorem.

We thus conclude that in the case 
of $h(0) = 0$ 
$$b = 1/\lim_{t\rightarrow 0} h^\prime(t).$$
with $b = 0$ if $h^\prime(t) \rightarrow \infty$ for $t \rightarrow 0$.
\\
$\mbox{ }$ \hfill $\Box$ \\

\section{Concluding remarks}

The proofs of equation~\eqref{duality} suggest that the relaxation modulus and the 
creep function should be considered as convolution operators. Such an approach allows 
incorporating the identity operator in this class. The identity operator cannot be ignored 
in the context of the duality equation \eqref{duality}. Even though one might set  $\beta = 0$ in \eqref{f0} 
a "Newtonian viscosity" term $b\, t$ can appear in equation~\eqref{ht}.

Another important operator in
viscoelastic theory is the integral $\mathrm{I}$
$$(\mathrm{I}\, g)(t) := \int_0^t g(s) \, \dd s$$
For example elastic media satisfy a constitutive equation of the form $\sigma = E \, 
\epsilon \equiv E\, \mathrm{I} \, \dot{\epsilon}$ for some $E > 0$. Unlike $\mathrm{U}$ the operator 
$\mathrm{I}$ can be expressed as a convolution with the unit step function, which is LICM.

Since $t\,\widetilde{} = 1/p^2$ and $(\mathrm{I}^2 \,f)(t) = \int_0^t f(s)\, \dd s$, 
equation \eqref{duality} can be expressed in an operator form
$$R \ast C \equiv C \ast R = \mathrm{I}^2$$
so that from the constitutive equation $\sigma = R \ast \dot{\epsilon}$ follows the inverse relation 
$C \ast \dot{\sigma} = C \ast \dd/\dd t [ R \ast \dot{\epsilon}]  = \epsilon$.

An anisotropic version of the duality relation eq.~\eqref{duality} was studied in 
\cite {HanDuality}. Our new method cannot be extended to matrix-valued functions because a theory of 
matrix-valued CBFs and Stieltjes functions has not been developed yet.

\appendix

\section{A few relevant properties of LICM, complete Bernstein and Stieltjes functions.}

For details, see \cite{SchillingAl,SerHanJMP}. 

In order to focus on those statements which 
are of use for us we shall consider as definitions some statements that appear as theorems  
in the references cited above.

An infinitely differentiable function $f$ is said to be LICM if it is completely monotone:
$$(-1)^n\, \D^n f(t) \geq 0 \     \text{for}  \      n = 0,1, 2,...$$
and integrable in a neighborhood of 0.

If $f(t)$ is LICM, then there is a real number $a \geq 0$ and a 
Borel measure $\mu$ on $]0,\infty[$ 
satisfying \eqref{x2} such that
$$f(t) = a + \int_{]0,\infty[} \e^{-r t}\, \mu(\dd r)$$ 
The last equation can also be 
recast in a more familiar form
$$f(t) = \int_{[0,\infty[} \e^{-r t}\, \mu(\dd r)$$
by defining $\mu(\{0\}) = a$. 

The inverse implication is also true.

A Bernstein function is defined as an infinitely differentiable function $g$ on $[0,\infty[$ satisfying the inequalities 
$g(t) \geq 0$ and $(-1)^n \,\D^n g(t) \leq 0$ for $n = 1,2,\ldots$ .

If $f$ is LICM, $a, b \geq 0$, then $g(t) = \int_0^t g(s) \, \dd s + a + b\, t$ 
is a Bernstein function. The derivative of a Bernstein function is LICM.

A function $f$ is a Stieltjes function if there are two real numbers $a, b \geq 0$ and
a Borel measure $\mu$ satisfying \eqref{x2} such that 
\begin{equation} \label{app1}
f(p) = a  + \frac{b}{p} + \int_{]0,\infty[}\frac{\mu(\dd r)}{p + r} 
\end{equation}
The integration extends over $0 < r < \infty$. The second term can be incorporated in the integral
by extending the integration to $[0,\infty[$, but we prefer to keep it separate.

A function $f$ is a CBF if there are two real numbers $a, b \geq 0$ and a Borel measure 
$\nu$ on $]0,\infty[$
satisfying \eqref{x0} 
such that 
\begin{equation}\label{app2}
f(p) = a + b\, p + p \int_{]0,\infty[} \frac{\nu(\dd r)}{p + r}
\end{equation}

$f(p)$ is a Stieltjes function if and only if $p\, f(p)$ is a CBF.
This statement follows from the integral representations of CBFs and Stieltjes functions above.

The following non-linear relation between the CBFs and the Stieltjes functions is the key to 
the proof above:
A function  $f(p)$ not identically zero  is a Stieltjes function if and only if  $1/f(p)$ is 
a CBF.

\end{document}